\documentclass[reqno,10pt]{article}
\usepackage{amsfonts}
\usepackage{hyperref}
\usepackage{amsmath}
\usepackage{amssymb}
\usepackage{mathrsfs}

\newtheorem{thm}{Theorem}[section]
\newtheorem{lemma}[thm]{Lemma}

\newtheorem{corol}[thm]{Corollary}
\newtheorem{propos}[thm]{Proposition}
\newtheorem{rema}{Remark}[section]
\def\bp{\begin{propos}}
\def\ep{\end{propos}}
\def\bt{\begin{thm}}
\def\et{\end{thm}}
\def\bco{\begin{corol}}
\def\eco{\end{corol}}
\def\bl{\begin{lemma}}
\def\el{\end{lemma}}
\def\br{\begin{rema}}
\def\er{\end{rema}}
\def\be{\begin{equation}}
\def\ee{\end{equation}}
\def\ba{\begin{array}}
\def\ea{\end{array}}
\def\bena{\begin{eqnarray}}
\def\eena{\end{eqnarray}}

\def\P{{\mathbb P}}
\def\E{{\mathbb E}}

\def\a{{\alpha}}
\def\d{{\delta}}\def\D{{\Delta}}

\def\QED{\hfill$\square$\vskip 3mm}

\def\Dp{\displaystyle}
\def\Df{\Dp\frac}
\def\hb{\hbox}
\def\({\left(}
\def\){\right)}

\begin{document}

\title{\LARGE On the Degree Sequence and its Critical Phenomenon of an Evolving Random Graph
Process\\[5mm]
 \footnotetext{Key
words and phrases: degree sequence; power law; critical phenomenon;
real-world networks.}}

\author{{Xian-Yuan Wu$^1$}\thanks{Supported in part
by the Natural Science Foundation of China},\ \ \ {Zhao
Dong$^2$}\thanks{Supported in part by the Natural Science Foundation
of China under grants 10671197 and 10721101},\ \ \ {Ke
Liu$^2$}\thanks{Partially supported by the Natural Science
Foundation of China under grants 60674082, 70221001 and 70731003.},\
\ \hb{ and }{Kai-Yuan Cai}$^3$\thanks{Supported by the Natural
Science Foundation of China and MicroSoft Research Asia under grant
60633010}} \vskip 10mm
\date{}
\maketitle {\small \vskip-20mm \begin{center}
\begin{minipage}{13cm} \noindent\hskip -2mm$^1$School of
Mathematical Sciences, Capital Normal University,
Beijing, 100037, China. Email: \texttt{wuxy@mail.cnu.edu.cn}\\[-5mm]

\noindent\hskip -2mm$^{2}$Academy of Mathematics and System
Sciences, Chinese Academy of Sciences, Beijing, 100190, China.
Email: \texttt{dzhao@amss.ac.cn;}\ \texttt{kliu@amss.ac.cn}\\[-5mm]

\noindent\hskip -2mm$^3$Department of Automatic Control, Beijing
University of Aeronautics and Astronautics, Beijing, 100083, China.
Email: \texttt{kycai@buaa.edu.cn}
\end{minipage}
\end{center}
\vskip 2mm
\begin{center} \begin{minipage}{14cm}
{\bf Abstract}: In this paper we focus on the problem of the degree
sequence for the following random graph process. At any time-step
$t$, one of the following three substeps is executed: with
probability $\a_1$, a new vertex $x_t$ and $m$ edges incident with
$x_t$ are added; or, with probability $\a-\a_1$, $m$ edges are
added; or finally, with probability $1-\a$, $m$ random edges are
deleted. Note that in any case edges are added in the manner of {\it
preferential attachment}. we prove that there exists a critical
point $\a_c$ satisfying: 1) if $\a_1<\a_c$, then the model has power
law degree sequence; 2) if $\a_1>\a_c$, then the model has
exponential degree sequence; and 3) if $\a_1=\a_c$, then the model
has a degree sequence lying between the above two cases.
\end{minipage}
\end{center}}

\vskip 5mm
\section{Introduction and statement of the results}
\renewcommand{\theequation}{1.\arabic{equation}}
\setcounter{equation}{0}

Recently there has been much interest in studying large-scale
real-world networks and attempting to model their properties. For a
general introduction to this topic, readers can refer to Albert and
Barab\'asi \cite{AB}, Aiello, Chung and Lu \cite{ACL}, Bollob\'as
and Riordan \cite{BR}, Hayes \cite{H}, Newman \cite{N} and Watts
\cite{W}. Although the study of real-world networks as graphs can be
traced back to long time ago such as the classical model proposed by
Erd\"os and R\'enyi \cite{ER} and Grilbert \cite{G}, recent
influential activity perhaps started with the work of Watts and
Strogatz about the `small-world phenomenon' published in 1998
\cite{WS}. Another influential work may be due to the scale-free
model proposed by Bollob\'as and Albert in 1999 \cite{BA}. Since
then various forms of scale-free phenomenon have been widely
revealed. In particular, power law degree distributions have been
extensively investigated. Many new models have been introduced to
circumvent the shortcomings of the classical models introduced by
Erd\"os and R\'enyi \cite{ER} and Grilbert \cite{G}. One class of
these new models was aimed to explain the underlying causes for the
emergence of power law degree distributions. This can be observed in
`LCD model' \cite{BR2} and its generalization due to Buckley and
Osthus \cite{BO}, `copying' models of Kumar {\it et al}.
\cite{KRRS}, the very general models defined by Copper and Frieze
\cite{CF} and the other model with random deletions defined by
Copper, Frieze and Vera \cite{CFV} {\it etc}.

For the real-world network of World Wide Web/Internet, experimental
studies by Albert, Barab\'asi and Jeong \cite{ABJ}, Broder {\it et
al}. \cite{BKM} and Faloutsos, Faloutsos and Faloutsos \cite{FFF}
demonstrated that the proportion of vertices of a given degree
follows an approximate inverse power law, i.e., the proportion of
vertices of degree $k$ is approximately $Ck^{-\a}$ for some
constants $C$ and $\a$. However other forms of the degree
distributions can also be observed in real-world networks (see
\cite{ASBS} and \cite{Str}). For example, Guassian distributions can
be observed in the acquaintance network of Mormons \cite{BKEMS};
exponential distribution can be observed in the powergrid of
southern California \cite{WS}. On the other hand, the degree
distribution of the network of world airports \cite{ASBS}
interpolates between Gaussian and exponential distributions, whereas
the degree distribution of the citation network in high energy
physics \cite{LLJ} interpolates between exponential and power law
distributions. For more forms of degree distributions, readers can
refer to~\cite{SAB}.

Different models often lead to different forms of degree
distributions. An interesting problem arises naturally: does it
exist some dynamically evolving random graph process which brings
forth various degree distributions by continuous changing of its
parameters only?
This phenomenon has been numerically investigated
in reference \cite{ZJW}: For a general model of collaboration
networks in~\cite{ZJW} ,
Zhou {\it et al.}  indicate that, while a relevant parameter $\a$
increases from $0$ to $1.5$, four kinds of degree distributions
appear as {\it exponential, arsy-varsy, semi-power law }and {\it
power law} in turn. Note that the above classification is rather
rough as no unambiguous borderline between two neighboring patterns
is determined. However, to the best of our knowledge, it seems that
the problem and its answer have not been formulated in a
mathematically rigorous manner. In this paper we focus on a model
with edge deletions and provide precise analysis, while a parameter
varies, the model exhibits various degree distributions.

Now, we begin to introduce our model and then state our main
results. Consider the following process which generates a sequence
of graphs $G_t=(V_t,E_t)$, $t\geq 1$. Write $v_t=|V_t|$ and
$e_t=|E_t|$.

\vskip 3mm {\it Time-Step 1.} {Let $G_1$ consist of an isolated
vertex $x_1$.}

\vskip 2mm {\it Time-Step $t\geq$ 2.} {

  1, With probability $\alpha_1>0$ we add a vertex $x_t$ to
  $G_{t-1}$. We then add $m$ random edges incident with $x_t$.
  In the case of $e_{t-1}>0$, the $m$ random neighbours
  $w_1,w_2,\ldots,w_m$ are chosen independently. For $1\leq i\leq m$
  and $w\in V_{t-1}$,
  \be\label{1.1}\P(w_i=w)=\frac{d_w(t-1)}{2e_{t-1}},\ee
  where $d_w(t-1)$ denotes the degree of vertex $w$ at the beginning
  of substep $t$. Thus neighbours are chosen by {\it preferential attachment.}
  In case of $e_{t-1}=0$,
  then we add a new vertex $x_t$ and join it to a randomly chosen
  vertex in $V_{t-1}$.

  2, With probability $\alpha-\alpha_1\geq 0$ we add $m$
  random edges to existing vertices. If $e_{t-1}>0$, then both endpoints are chosen
  independently with the same probabilities as in (\ref{1.1}). Otherwise, we do nothing.

  3, With probability $1-\alpha\geq 0$ we delete $\min\{m,e_{t-1}\}$ randomly chosen edges from $E_{t-1}$.}

\vskip 5mm

\br\label{r0} The deference between our model and the model
introduced in \cite{CFV} is that, in our setting, vertex deletions,
loop and multi-edge erasures are forbidden, which makes $\{e_t:t\geq
1\}$ Markovian and makes it possible for us to give exact estimation
to $e_t$. \er

In order to make the problem meaningful, the following inequalities
are natural and necessary: \bena && 1/2<\a\leq 1; \ \ 0<\a_1\leq
\a.\label{1.3}\eena For given $\a$ and $\a_1$ satisfying
(\ref{1.3}), define \be\label{1.6'}\a_c:=4\a-2,\ \ \eta:=\a_cm/2,\ee
and choose $\epsilon=\epsilon(\a,\a_1)\in(0,\eta)$ such that
\be\label{2.3}\rho_\epsilon:=\max\left\{\frac{m(\a_c-\a_1)}{2(\eta-\epsilon)},\frac
12\right\}<1.\ee Note that in case of $\a_1\geq \a_c$,
$\rho_\epsilon=\frac 12$. Let
\be\label{1.6}\beta=\Dp\frac{\a_c}{\a_c-\a_1};\ \
\gamma=1-\Df{\a_1-\a_c}{2(1-\a)};\ \
\theta=\frac{2\a_c-\a_1}{2\a_c};\ \ \mu=\frac{\a_c}{2(1-\a)}.\ee
Obviously, $\beta$ is well defined when $\a_1\not=\a_c$ and
$0<\gamma<1$ when $\a_1>\a_c$. To get our main results, besides
(\ref{1.3}), the following condition is necessary
\be\label{1.4}\a_1<2\a_c.\ee

Now, Let $D_k(t)$ be the number of vertices with degree $k\geq 0$ in
$G_t$ and let $\overline{D}_k(t)$ be the expectation of $D_k(t)$.
The main results of this paper follow as

\bt\label{th1} Assume that (\ref{1.3}) and (\ref{1.4}) hold. Then
$\a_c$ defined in (\ref{1.6'}) is a critical point for the degree
sequence of the model satisfying:

1) if $\a_1<\a_c$, then there exists a constant $C_1=C_1(m,\a,
\a_1)$ such that, for any $\nu\in(0,1-\rho_\epsilon)$,
\be\label{1.7}\left|\Df{\overline{D}_k(t)}{t}-C_1k^{-1-\beta}\right|=
O(t^{\rho_\epsilon+\nu-1})+O(k^{-2-\beta});\ee

 2) if $\a_1>\a_c$, then there exists a constant
$C_2=C_2(m,\a,\a_1)$ such that
\be\label{1.8}\left|\Df{\overline{D}_k(t)}{t}-C_2\gamma^k
k^{-1+\beta}\right|=O(t^{-\theta})+O(\gamma^kk^{-2+\beta});\ee

 3) if $\a_1=\a_c$, then
there exists a constant $C_c=C_c(m,\a,\a_1)$ such that, for any
$\nu\in(0,\frac12)$,
\be\label{1.9}\left|\Df{\overline{D}_k(t)}{t}-C_cu_c(k)\right|=O
(t^{-\frac 12+\nu})\ee uniformly in $k$.

Where $u_c(k)=\Dp\int^1_0\Dp t^{k-1}\Dp e^{-\frac\mu{1-t}}dt$ and
$\beta$, $\gamma$, $\theta$ and $\mu$ are given in (\ref{1.6}). \et

\br\label{r1} The integral $u_c(k)=\Dp\int^1_0\Dp t^{k-1}\Dp
e^{-\frac\mu{1-t}}dt$ can be rewritten as \bena
u_c(k)&\hskip-3mm=&\hskip-3mm\left[\sum_{i=0}^{k-2}\sum_{l=0}^{k-2-i}
\binom {k-1}i\Df{(k-i-l-1)!}{(k-i)!}(-1)^{k-i-l-1}\right]e^{-\mu}\nonumber\\[3mm]
&&\hskip-3mm+\left[\sum_{i=0}^{k-1}\binom{k-1}i\Df{\mu^{k-i-1}}{(k-i)!}\right]\int_1^{+\infty}t^{-2}e^{-\mu
t}dt.\nonumber\eena With help of computer calculation, $u_c(k)$
satisfies $$\lim_{k\rightarrow \infty}{\ln
u_c(k)}/{(-k)}=\lim_{k\rightarrow \infty}{(-\ln k)}/{\ln
u_c(k)}=0.$$\er

Based on Theorem \ref{th1}, we can obtain following two corollaries,
which provide a complete distinction with respect to the parameters
between the degree sequences for the present model.

\bco\label{c1}If the parameters satisfy that

1) $\a> 2/3$; or

2) $\a\leq 2/3$ and $\a_1<\a_c$, \\
then the present random graph process has the power law degree
sequence (\ref{1.7}). \eco

\bco\label{c2}Assume $\a\leq 2/3$.

1) If $\a_c<\a_1<2\a_c$, then the present random graph process has
the exponential degree sequence  (\ref{1.8}).

2) If $\a_1=\a_c$, then the present random graph process has the
critical degree sequence (\ref{1.9}). \eco

\br\label{r2} When $\a>2/3$, for any $\a_1$,  the inequality
$\a_1\leq \a<\a_c=4\a-2$ holds always, therefore, the part 1) of
Corollary \ref{c1} follows from the part 1) of Theorem \ref{th1}.
The part 2) of Corollary \ref{c1} and Corollary \ref{c2} are
straightforward from Theorem \ref{th1}.\er

\br\label{r3} A special case of the part 1) in Corollary \ref{c1} is
$\a=\a_1=1$. In this case, the model has a power law degree sequence
as $Ck^{-3}$, which coincides with the result of \cite{BRST}.
Furthermore, for any $\a \in (1/2, 1]$ and $\a_1=2\a-1$, the model
has the degree sequence $Ck^{-3}$.\er

\br\label{r4}The results are unclear for the following case: $\a\leq
2/3$, $2\a_c\leq \a_1\leq\a$. Clearly, this case can only appear
when $\a\leq 4/7$. It is natural to conjecture that the model
possesses an exponential degree sequence in this case. \er

The methodology of the proof for the main results follows the
standard procedure which can be found in \cite{CF} and \cite{CFV}.
The rest of the paper is organized as follows. In Section 2, we
bound the degree of vertex in $G_t$. In Section 3, we establish the
recurrence for $\overline{D}_k(t)$ and then derive the approximation
of $\overline{D}_k(t)$ by a recurrence with respect to $k$. Finally,
in section 4, we solve the recurrence in $k$ using Laplace's method
\cite{J} and finish the proof of Theorem~\ref{th1}.

\section{Bounding the Degree}
\renewcommand{\theequation}{2.\arabic{equation}}
\setcounter{equation}{0}

For times $s$ and $t$ with $1\leq s\leq t$, let $d_{x_s}(t)$ be the
degree of vertex $x_s$ in $G_t$. If $x_s$ is not added in Time-Step
$s$, i.e., at Time-Step $s$, one of the other two substeps is
executed, put $d_{x_s}(t)=0$. In this section, we will concentrate
on the upper bound of $d_{x_s}(t)$.

For the present model, the estimation for $v_t$ is derived in
\cite{CFV} as \bena |v_t-\a_1 t|\hskip -4mm&&\leq ct^{1/2}\log t,\ \
\ \mbox{qs},\nonumber\eena for any constant $c>0$. We say an event
happens {\it quite surely} (qs) if the probability of the
complimentary set of the event is $O(t^{-K})$ for any $K>0$.

For the estimation of $e_t$, it can be derived by the same argument
as in \cite{CFV} that \bena|e_t-\eta t|\hskip -4mm&&\leq
ct^{1/2}\log t,\ \ \ \mbox{qs},\label{2.1}\eena for any constant
$c>0$.

By a standard argument on large deviation (see e.g. \cite{LS} and
\cite{S}), one further has: for any $\epsilon
>0$, there exists $c_1, c_2>0$ such that \be\label{2.2'}\P(e_t\leq
(\eta-\epsilon)t)\leq c_1\exp\{-c_2t\},\ee for all $t\geq 1$.

The following is our bounding for $d_{x_s}(t)$, note that our result
is based on the exact estimation (\ref{2.2'}) for $e_t$. In our
opinion, to bound the degree of vertex effectively, aforehand good
estimations for $e_t$ are necessary.

\bl\label{l1} For any $\a\in(1/2,1]$ and $\a_1\in (0,\a]$,
\be\label{2.2}d_{x_s}(t)\leq(t/s)^{\rho_\epsilon}(\log t)^3 \ \ \
\mbox{qs},\ee\el where $\rho_\epsilon$ is given in (\ref{2.3}).

{\it Proof}: Fix $s\leq t$, suppose that $x_s$ is added in Time-Step
$s$. Let $X_\tau=d_{x_s}(\tau)$ for $\tau=s,s+1,\ldots,t$ and let
\be\label{2.0}\lambda=\frac{(s/t)^{\rho_\epsilon}}{NM_\epsilon(\log
t+1)}, \ee where $N$ be large enough and will be determined later,
and $M_\epsilon=\Df{12m^2}{\eta-\epsilon}$. Let $Y$ be the
$\{1,2,3\}$-valued random variable with $\P(Y=1)=\a_1$,
$\P(Y=2)=\a-\a_1$ and $\P(Y=3)=1-\a$. Then conditional on $X_\tau=x$
and $e_\tau\geq m$, we have \be\label{2.4}X_{\tau+1}=
x+I_{\{Y=1\}}B(m,\frac x{2e_\tau})+I_{\{Y=2\}}B(2m,\frac
x{2e_\tau})-I_{\{Y=3\}}S(m,\frac x{e_\tau}),\ee where $B(m,p)$ is
the Binomial random variable with parameter $(m,p)$ and
$S(m,\Dp\frac x{e_\tau})$ is the super geometric random variable
with parameter $(e_\tau,x,m)$.

Noticing that $\lambda$ is small enough for large $N$, using the
basic inequality
\[e^{-y}\leq 1-y+2y^2 \hb{ for small } y>0\] and the fact that
$S(m,\Df{x}{e_\tau})\leq m$, (\ref{2.4}) implies
\bena &&\E\(\E(e^{\lambda X_{\tau+1}}\mid X_\tau=x,e_\tau)\mid e_\tau\geq m\)\nonumber\\[3mm]&&\hskip
-5mm\leq e^{\lambda
x}\left\{\a_1\left[1+\frac{x}{2e_\tau}(e^{\lambda}-1)\right]^m+(\a-\a_1)
\left[1+\frac{x}{2e_\tau}(e^\lambda-1)\right]^{2m}\right.\nonumber\\[3mm]&&\hskip9mm\left.+(1-\a)\left[1-\lambda\E\(
S(m,\frac{x}{e_\tau})\)+2\lambda^2m\E\(S(m,\frac{x}{e_\tau})\)\right]\right\}\label{2.5}.\eena
Using the inequalities \[e^y\leq 1+y+2y^2 \hb{ for small $y>0$}\]
and
\[(1+y)^m\leq 1+my+\frac{m^2}2y^2 \hb{ for small }y>0\] to the right
hand side of (\ref{2.5}) in turn, we get \bena &&\E\(\E(e^{\lambda
X_{\tau+1}}\mid X_\tau=x,e_\tau)\mid e_\tau\geq m\)\leq e^{\lambda
x}\left\{1+\frac{m\lambda
x}{2e_\tau}\left[(\a_c-\a_1)+12m\lambda\right]\right\}\nonumber\\[3mm]
&&\hskip10mm\leq e^{\lambda x}\left\{1+\frac{m\lambda
x}{2e_\tau}\left[\max\left\{(\a_c-\a_1),\frac{\eta-\epsilon}m\right\}+12m\lambda\right]\right\}\label{2.5'}\\[3mm]
&&\hskip10mm\leq e^{\lambda x}\left\{1+\frac{m\lambda
x\max\left\{(\a_c-\a_1),\frac{\eta-\epsilon}m\right\}}{2e_\tau}(1+M_\epsilon\lambda)\right\}\nonumber\\[3mm]
&&\hskip10mm\leq \exp\left\{\lambda
x\left[1+\frac{\max\left\{(\a_c-\a_1),\frac{\eta-\epsilon}m\right\}m}{2e_\tau}(1+M_\epsilon\lambda)\right]\right\}.\label{2.6}\eena

Now, we express $\E(e^{\lambda X_{\tau+1}}\mid X_\tau=x)$ as \bena
&&\hskip 4mm\E\(e^{\lambda X_{\tau+1}}\mid
X_\tau=x\)=\E\(\E(e^{\lambda X_{\tau+1}}\mid
X_\tau=x,e_\tau)\)\nonumber\\[3mm]
&&=\E\(\E(e^{\lambda X_{\tau+1}}\mid X_\tau=x,e_\tau)\mid
e_\tau<m\)\P(e_\tau<m)\nonumber\\[3mm]&&\hskip 4mm+\E\(\E(e^{\lambda X_{\tau+1}}\mid X_\tau=x,e_\tau)\mid e_\tau\geq
m\)\P(e_\tau\geq m)\nonumber\\[3mm]
&&=:I+II.\label{2.6'}\eena On one hand, conditional on $X_\tau=x$
and $e_\tau<m$, $X_{\tau+1}\leq x+m\leq e_\tau+m\leq 2m$ holds
always, so \be\label{2.7}I\leq e^{2m}\P(e_\tau<m).\ee On the other
hand, $II$ can be expressed as \bena &&II=\E\left(\E(e^{\lambda
X_{\tau+1}}\mid X_\tau=x,e_\tau)\mid e_\tau\geq m,e_\tau\geq
(\eta-\epsilon)\tau\right)\nonumber\\[2mm]&&\hskip
8mm\times\P\left(e_\tau\geq m,e_\tau\geq
(\eta-\epsilon)\tau\right)\nonumber\\[2mm]
&&\hskip 8mm+\E\left(\E(e^{\lambda X_{\tau+1}}\mid
X_\tau=x,e_\tau)\mid
e_\tau\geq m,e_\tau<(\eta-\epsilon)\tau\right)\nonumber\\[2mm]&&\hskip
8mm\times\P\left(e_\tau\geq
m,e_\tau<(\eta-\epsilon)\tau\right),\nonumber\eena by (\ref{2.6})
and the fact that $x\leq e_\tau$, \bena &II&\leq \exp\left\{\lambda
x\left[1+\frac{\max\{(\a_c-\a_1),\frac{\eta-\epsilon}m\}m}{2(\eta-\epsilon)\tau}(1+M_\epsilon\lambda)\right]\right\}\nonumber\\[3mm]
&&\hskip 3mm
+e^{\lambda(\eta-\epsilon)\tau}\exp\left\{\frac{\max\{(\a_c-\a_1),\frac{\eta-\epsilon}m\}m\lambda}{2}(1+M_\epsilon\lambda)\right\}
\P(e_\tau<(\eta-\epsilon)\tau)\nonumber\\[3mm]
&&\leq\exp\left\{\lambda
x\left[1+\frac{\rho_\epsilon}{\tau}(1+M_\epsilon\lambda)\right]\right\}+C'
e^{\lambda(\eta-\epsilon)\tau}\P(e_\tau<(\eta-\epsilon)\tau)\eena
for some constant $C'=C'(\a,\a_1,\epsilon,m)>0$. By (\ref{2.2'}) and
(\ref{2.0}), choosing $N$ large enough, then there exists constants
$c_3,c_4>0$ such that \be\label{2.8} II\leq \exp\left\{\lambda
x\left[1+\frac{\rho_\epsilon}{\tau}(1+M_\epsilon\lambda)\right]\right\}+c_3\exp\{-c_4\tau\}.\ee

Combining (\ref{2.6'})-(\ref{2.8}), using (\ref{2.2'}) again for
(\ref{2.7}), then there exists constants $c_5,c_6>0$ such that
\be\label{2.9}\E(e^{\lambda X_{\tau+1}}\mid X_\tau=x)\leq
\exp\left\{\lambda
x\left[1+\frac{\rho_\epsilon}{\tau}(1+M_\epsilon\lambda)\right]\right\}+c_5\exp\{-c_6\tau\}.\ee
Thus \be\label{2.10}\E(e^{\lambda X_{\tau+1}})\leq
\E\left(\exp\left\{X_\tau\lambda\left(1+\frac{\rho_\epsilon(1+M_\epsilon
\lambda)}{\tau}\right)\right\}\right)+c_5\exp\{-c_6\tau\}.\ee

Now, put $\lambda_t=\lambda$ and
$\lambda_{\tau-1}=\lambda_\tau(1+\frac{\rho_\epsilon(1+M_\epsilon
\lambda_\tau)}{\tau})$. Obviously, if $\lambda_s$ is small enough,
then (\ref{2.10}) holds for $\lambda_{\tau+1},
\tau=s,s+1,\ldots,t-1$. This will imply that
\be\label{2.11}\E(e^{\lambda X_t})=\E(e^{\lambda_t X_t})\leq
\E(e^{\lambda_sX_s})+c_5\sum_{\tau=s}^t\exp\{-c_6\tau\}\leq
e^{m\lambda_s}+C''\ee for some constant $C''>0$.

Let $\Lambda=\Dp\frac{10}{NM_\epsilon(\log t+1)}$, note that
$\Lambda$ can be taken small enough uniformly in $t$ by taking $N$
large enough. Now provided $\lambda_\tau\leq {\Lambda}$, we can
write
\[\lambda_{\tau-1}\leq \lambda_\tau\left(1+\frac{\rho_\epsilon(1+M_\epsilon
\Lambda)}{\tau}\right)\] and then
\bena\lambda_s\hskip-5mm&&\leq\lambda\prod_{\tau=s}^t\left(1+\frac{\rho_\epsilon(1+M_\epsilon
\Lambda)}{\tau}\right)\leq 10
\lambda(t/s)^{\rho_\epsilon}\nonumber\eena which is $\leq\Lambda$ by
the definition of $\lambda$.

Put $u=(t/s)^{\rho_\epsilon}(\log t)^3$, by (\ref{2.11}) we get
\[\P(X_t\geq u)\leq (e^{m\lambda_s}+C'')e^{-\lambda u}=O(t^{-K})\nonumber\] for any constant $K>0$ and the
Lemma follows. \QED

\br For any $n$ large enough, $\rho_\epsilon$ can be retaken as
$\Dp\max\left\{\frac{m(\a_c-\a_1)}{2(\eta-\epsilon)},\frac
1n\right\}$, in fact, this can be done by enlarging $\a_c-\a_1$ to
$\Dp\max\left\{(\a_c-\a_1), \frac{2(\eta-\epsilon)}{nm}\right\}$
instead in (\ref{2.5'}). Thus, in the case of $\a_1\geq \a_c$,
$\rho_\epsilon$ can be taken as $1/n$. Certainly, if this is done as
above, constants $M_\epsilon$ and $N$ should be retaken
correspondingly.\er

\section{The recurrence for $\overline D_k(t)$}
\renewcommand{\theequation}{3.\arabic{equation}}
\setcounter{equation}{0}

In this Section, we follow the basic procedures in \cite{CFV} to
establish the recurrence for $\overline D_k(t)$. Put $D_{-1}(t)=0$
for all $t\geq 1$. For $k\geq 0$, we have \bena
&&\overline D_k(t+1)=\overline D_k(t)\nonumber\\[3mm]&&\hskip
2mm+(2\a-\a_1)m\E\left(\left.-\frac{kD_k(t)}{2e_t}+\frac{(k-1)D_{k-1}(t)}{2e_t}+O\(\Df{\D_t}{e_t}\)\right|
e_t>0\right)\P(e_t>0)\nonumber\\[3mm]
&&\hskip
2mm+(1-\a)m\E\left(\left.\frac{(k+1)D_{k+1}(t)}{e_t}-\frac{kD_k(t)}{e_t}+O\(\Df{\D_t}{e_t}\)\right|
e_t\geq m\right)\P(e_t\geq
m)\nonumber\\[3mm]
&&\hskip
2mm+\a_1I_{k=m}\P(e_t>0)+O(\P(e_t=0))+O(\P(e_t<m)).\label{3.0}\eena
Here $\D_t$ denotes the maximum degree in $G_t$ and the term
$O\(\Df{\D_t}{e_t}\)$  accounts for the probability that we create
larger than one degree changes for some vertices at Time-Step $t+1$.
By (\ref{2.2'}) and Lemma \ref{l1}, we have \be\label{3.0'}
\Df{\D_t}{e_t}\leq O(t^{\rho_\epsilon-1}(\log t)^3),\ \ \ qs.\ee

The term $\Dp\E\(\left.\frac{kD_k(t)}{e_t}\right|e_t>0\)$ can be
expressed as \bena&&\E\left(\left.\frac{kD_k(t)}{e_t}\right|
e_t>0\right)\nonumber\\[3mm]&&=\E\(\left.\frac{kD_k(t)}{e_t}\right| |e_t-\eta
t|\leq
t^{1/2}\log t\)\P\(|e_t-\eta t|\leq t^{1/2}\log t\mid
e_t>0\)\nonumber\\[3mm]
&&\hskip 4mm+\E\(\left.\frac{kD_k(t)}{e_t}\right| |e_t-\eta t|>
t^{1/2}\log t,e_t>0\)\P(|e_t-\eta t|> t^{1/2}\log t\mid
e_t>0)\nonumber\\[3mm]
&&=\frac{\E\({kD_k(t)}\mid |e_t-\eta t|\leq t^{1/2}\log
t\)\P(|e_t-\eta t|\leq t^{1/2}\log t\mid e_t>0)}{\eta
t}\nonumber\\[2mm]&&\hskip 5mm\times(1+O(t^{-1/2}\log t))+O(\P(|e_t-\eta
t|> t^{1/2}\log t\mid e_t>0)),\label{3.1'}\eena where we used  the
fact that $kD_k(t)\leq 2e_t$ to hand the second term. For $k\geq 1$,
we have $\overline D_k(t)=\E(D_k(t)\mid e_t>0)\P(e_t>0)$, so
\bena&&\E({kD_k(t)}\mid |e_t-\eta t|\leq t^{1/2}\log t)\P(|e_t-\eta
t|\leq t^{1/2}\log t\mid
e_t>0)\nonumber\\[2mm]&&=k\overline D_k(t)-\E({kD_k(t)}\mid |e_t-\eta t|>
t^{1/2}\log t,e_t>0)\nonumber\\[2mm]&&\hskip 5mm\times\P(|e_t-\eta t|>
t^{1/2}\log t\mid e_t>0)\nonumber\\[2mm]&&=k\overline
D_k(t)+O(t\cdot\P(|e_t-\eta t|> t^{1/2}\log t\mid
e_t>0)).\label{3.2'}\eena Thus, using (\ref{2.1}), we have for
$k\geq 0$ \be\label{3.3'}\E\(\left.\frac{kD_k(t)}{e_t}\right|
e_t>0\)=\frac{k\overline D_k(t)}{\eta t}+O(t^{-1/2}\log t).\ee
Similarly, \be\label{3.4''}\E\(\left.\frac{kD_k(t)}{e_t}\right|
e_t\geq m\)=\frac{k\overline D_k(t)}{\eta t}+O(t^{-1/2}\log t).\ee

Substituting (\ref{3.0'}), (\ref{3.3'}) and (\ref{3.4''}) into
(\ref{3.0}), using (\ref{2.2'}) again to the other terms, we derive
the following approximate recurrence for $\overline D_k(t)$:
$\overline D_{-1}(t)=0$ for all $t>0$ and for $k\geq 0$
\bena\overline D_k(t+1)&&\hskip-5mm=\overline
D_k(t)+(A_2(k+1)+B_2)\Df{\overline
D_{k+1}(t)}{t}+(A_1k+B_1+1)\Df{\overline D_{k}(t)}{t}\nonumber\\[-1mm]&&
\nonumber  \\[-1mm]
&&\hskip-2mm+(A_0(k-1)+B_0)\Df{\overline
D_{k-1}(t)}{t}+\a_1I_{k=m}+O(t^{\rho_\epsilon-1}(\log t)^3),
\label{3.1} \eena where
\[ A_2=\Df{1-\a}{2\a-1};\ A_1=-\Df{2-\a_1}{2(2\a-1)};\
A_0=\Df{2\a-\a_1}{2(2\a-1)};\ B_2=B_0=0\hb{ and }B_1=-1.\] Note that
the hidden constant, write as $L$, in term
$O(t^{\rho_\epsilon-1}(\log t)^3)$ of (\ref{3.1}) is uniform in $k$,
which follows from the fact that $e_t=O(t)$ and $k D_k(t)\leq
2e_t=O(t)$ uniformly in $k$.

If we heuristically put $\bar d_k=\Df{\overline D_k(t)}{t}$ and
assume it is a constant, we get \bena\bar
d_k&&\hskip-5mm=(A_2(k+1)+B_2)\bar
d_{k+1}+(A_1k+B_1+1)\bar d_{k}\nonumber\\[2mm]
&&\hskip-2mm+(A_0(k-1)+B_0)\bar
d_{k-1}+\a_1I_{k=m}+O(t^{\rho_\epsilon-1}(\log t)^3).\nonumber \eena
This leads to the consideration of the recurrence in $k$: $d_{-1}=0$
and for $k\geq -1$, \be\label{3.3}(A_2(k+2)+B_2)
d_{k+2}+(A_1(k+1)+B_1)d_{k+1} +(A_0k+B_0)d_{k}=-\a_1I_{k=m-1}.\ee

The following Lemma shows that, on certain conditions, (\ref{3.3})
is a good approximation to (\ref{3.1}). Note that our Lemma is a
generalization of Lemma~5.1 in \cite{CFV}.

\bl\label{l3} Let $d_k$ be a solution for (\ref{3.3}) such that
$|d_k|\leq \frac Ck$ for $k>0$ and a constant $C$. We have

1) if $\a_1\leq \a_c$, then, for any $\nu\in (0,1-\rho_\epsilon)$,
there exists a constant $M_1>0$ such that \be\label{3.4}|\overline
D_k(t)-td_k|\leq M_1t^{\rho_\epsilon+\nu},\ee for all $t\geq 1$ and
$k\geq-1$;

2) if $\a_c<\a_1<2\a_c$, then there exists a constant $M_2>0$ such
that \be\label{3.4'}|\overline D_k(t)-td_k|\leq M_2t^{1-\theta},\ee
for all $t\geq 1$ and $k\geq-1$, where $\theta$ is given in
(\ref{1.6}).\el

{\it Proof.} Let $\Theta_k(t)=\overline D_k(t)-td_k$ and
$k_0=k_0(t)=\lfloor t^{\rho_\epsilon}(\log t)^3\rfloor$. Lemma
\ref{l1} implies \be\label{3.5}0\leq\overline D_k(t)\leq t^{-10}\ \
\makebox{for }k\geq k_0(t).\ee

{\it Proof of part 1)}: Equation (\ref{3.5}) and $d_k\leq C/k$ imply
that (\ref{3.4}) holds for $k\geq k_0$ uniformly, i.e., there exists
a constant $N_1>0$, independent to $k$ and $t$, such that
$$|\overline D_k(t)-td_k|=|\Theta_k(t)|\leq N_1t^{\rho_\epsilon}$$ for all $k\geq
k_0(t)$ and $t\geq 1$.

Recall that the hidden constant in $O(t^{\rho_\epsilon-1}(\log
t)^3)$ of (\ref{3.1}) is denoted by $L$. For any
$\nu\in(0,1-\rho_\epsilon)$, let $R\geq L$ satisfying
$$Lt^{\rho_\epsilon-1}(\log t)^3\leq R t^{\rho_\epsilon+\nu-1}$$ for all $t\geq 1$. Let
$N_2=\frac R{\rho_\epsilon+\nu}+1$, take $\sigma>0$ such that
\be\label{31}1-\frac {R}{N_2}-(1+\sigma)\(1-\rho_\epsilon-\nu\)\geq
0,\ee and take $\d\in(0,1)$ such that
\be\label{32}\d^{1+\sigma}<e^{-1}<\d.\ee Let $t_1>0$ be an integer
such that
\be\label{32'}k_0(t)\leq\frac{-1}{A_1}t=\frac{2(2\a-1)}{2-\a_1}t\ee
and
\be\label{33}\d^{1+\sigma}\leq\left(1-\frac{1}{t+1}\right)^{{t+1}}
,\ \ \ \left(1-\frac{1-R/l}{t+1}\right)^{\frac{t+1}{1-R/l}}\leq \d
\ee for all $t\geq t_1$ and $l\geq N_2$.

Now, for the above $t_1$, let $N_3\ge N_1$ satisfying
 \be\label{3.6}|\Theta_k(t)|\leq N_3t^{\rho_\epsilon+\nu}\ \
\makebox{for all } 1\leq t\leq t_1 \makebox{ and } k\geq -1.\ee Take
\be\label{3.7} M_1=\max\{N_2, N_3\}.\ee We will prove that
(\ref{3.4}) holds for the above $M_1$ by induction. Our inductive
hypothesis is \[{\cal H}^1_t:\ \  |\Theta_k(t)|\leq
M_1t^{\rho_\epsilon+\nu}\ \ \hb{for all}\ \ k\ge -1.\] Note that
(\ref{3.6}) and (\ref{3.7}) imply that ${\cal H}^1_t$ holds for
$1\leq t\leq t_1$.

It follows from (\ref{3.1}) and (\ref{3.3}) that
\bena\Theta_k(t+1)&&\hskip-5mm=\Theta_k(t)+A_2(k+1)\Df{\Theta_{k+1}(t)}t+(A_1k+B_1+1)\Df{\Theta_k(t)}t\nonumber\\[-2mm]
&&\nonumber \\[-2mm]
&&\hskip-5mm+A_0(k-1)\Df{\Theta_{k-1}}{t}+O(t^{\rho_\epsilon-1}(\log
t)^3). \label{3.8}\eena For $t\geq t_1$, by (\ref{32'}), we have
$t+A_1k+B_1+1\geq 0$ and then (\ref{3.8}) implies
\bena|\Theta_k(t+1)|&&\hskip-5mm\leq
A_2(k+1)\Df{|\Theta_{k+1}(t)|}{t}+(t+A_1k+B_1+1)\Df{|\Theta_k(t)|}t\nonumber\\[2mm]
&&+A_0(k-1)\Df{|\Theta_{k-1}(t)|}{t}+Rt^{\rho_\epsilon+\nu-1}\nonumber\\[1mm]
&&\hskip-5mm\leq(t+A_2(k+1)+A_1k+B_1+1+A_0(k-1))M_1t^{\rho_\epsilon+\nu-1}+Rt^{\rho_\epsilon+\nu-1}\nonumber\\[2mm]
&&\hskip-5mm=(t+A_2+B_1+1-A_0)M_1t^{\rho_\epsilon+\nu-1}+Rt^{\rho_\epsilon+\nu-1}.\nonumber
\eena Let $\varepsilon_0=A_2+B_1+1-A_0=(\a_1-\a_c)/\a_c$, noticing
that $\a_1\leq \a_c$, we have $\varepsilon_0\leq 0$. Then, combining
(\ref{31}), (\ref{33}) and (\ref{3.7}), we have \bena
&&\Df{(t+\varepsilon_0)M_1t^{\rho_\epsilon+\nu-1}+Rt^{\rho_\epsilon+\nu-1}}{M_1(t+1)^{\rho_\epsilon+\nu}}\leq
\Df{M_1t^{\rho_\epsilon+\nu}+Rt^{\rho_\epsilon+\nu-1}}{M_1(t+1)^{\rho_\epsilon+\nu}}\nonumber\\[2mm]
&&=\left\{\left(1-\frac{1-
R/{M_1}}{t+1}\right)^{\frac{t+1}{1-R/M_1}}\right\}^{\frac{1-R/M_1}{t+1}}
\left\{\left(1-\frac{1}{t+1}\right)^{{t+1}}\right\}^{\frac{-(1-\rho_\epsilon-\nu)}{t+1}}\nonumber\\[4mm]
&&\Dp\leq
\d^{\frac{1-R/M_1}{t+1}}\cdot(\d^{1+\sigma})^{\frac{-(1-\rho_\epsilon-\nu)}{t+1}}=\d^{(1-\frac
R{M_1}-(1+\sigma)(1-\rho_\epsilon-\nu))/(t+1)}\nonumber\\[3mm]
&&\leq 1.\nonumber\eena The induction hypothesis ${\cal H}^1_{t+1}$
has been verified and the proof of part 1) is completed.

{\it Proof of part 2)}: In this case, we have $\a_c<\a_1<2\a_c$ and
then, for some $\nu\in(0,1/2)$, $\varepsilon_0\leq
\rho_\epsilon+\nu<1-\theta$ ( note that in this case
$\rho_\epsilon=1/2$). Same as what we have done for part 1), for
certain $\sigma>0$ and $\d\in(e^{-1},1)$, we have
\bena&&\Dp\Df{(t+\varepsilon_0)M_2t^{-\theta}+Rt^{-\theta}}{M_2(t+1)^{1-\theta}}
\leq\Dp\d^{{(1-\varepsilon_0-\frac
R{M_2}-(1+\sigma)\theta)}/{(t+1)}}\leq1\nonumber\eena for sufficient
large $t$ and $M_2$. This is enough for a inductive proof of
(\ref{3.4'}). \QED

\br\label{r5}$^{\makebox{\rm [Remark 5.2 in \cite{CFV}]}}$ Lemma
\ref{l3} implies that if there is a solution for (\ref{3.3}) such
that $d_k\leq C/k$, then $\lim_{t\rightarrow\infty}\overline
D_k(t)/t$ exists and equals to $d_k$. In particular, it is shown
that: if there exists a solution for (\ref{3.3}) such that $d_k\leq
C/k$, then the solution is unique.\er

\section{Solving (\ref{3.3}) and the proof of Theorem \ref{th1}}
\renewcommand{\theequation}{4.\arabic{equation}}
\setcounter{equation}{0}

In order to solve (\ref{3.3}), let us consider the following
homogeneous equation \be\label{4.1}(A_2(k+2)+B_2)
f_{k+2}+(A_1(k+1)+B_1)f_{k+1} +(A_0k+B_0)f_{k}=0,\ \ k\geq 1\ee
which is solved by Laplace's method as explained in \cite{J}.

For $k\geq 1$, we construct function $f_k$ has the following form
\be\label{4.2}f_k=\int_{a}^{b}t^{k-1}v(t)dt,\ee where constants $a$
and $b$, and function $v(t)$ are to be determined later.

Integrating by parts
\be\label{4.3}kf_k=[t^kv(t)]^b_a-\int_{a}^{b}t^{k}v'(t)dt.\ee Let
$$\phi_1(t)=A_2t^2+A_1t+A_0,\ \ \phi_0(t)=B_2t^2+B_1t+B_0.$$
Substituting (\ref{4.2}) and (\ref{4.3}) into (\ref{4.1}), we obtain
\be\label{4.4}[t^k\phi_1(t)v(t)]^b_a-\int_a^bt^k\phi_1(t)v'(t)dt+\int_a^bt^{k-1}\phi_0(t)v(t)dt=0.\ee
Equation (\ref{4.1}) will be satisfied if we have
\be\label{4.5}\Df{v'(t)}{v(t)}=\Df{\phi_0(t)}{t\phi_1(t)},\ee and
\be\label{4.6}[t^kv(t)\phi_1(t)]^b_a=0.\ee Let $a=0$ and $b$ equal
to a root of $v(t)\phi_1(t)=0$, the parameters $a$ and $b$ can be
determined satisfying (\ref{4.6}).

Obviously, $\phi_0(t)$ and $\phi_1(t)$ can be rewritten as
\be\label{4.7} \phi_0(t)=-t;\ \
\phi_1(t)=At^2-(A+B)t+B=A(t-1)(t-B/A),\ee where\bena
&&A=\Df{1-\a}{2\a-1},\ \ \
B=\Df{2\a-\a_1}{\a_c}=A+\frac{\a_c-\a_1}{\a_c}.\label{4.8}\eena

Now, we solve the equation (\ref{4.1}) in the following cases: 1),
$\a_1<\a_c$; 2), $\a_1>\a_c$ and 3), $\a_1=\a_c$ respectively.

For case $\a_1<\a_c$, we have $B>A$, then the differential equation
(\ref{4.5}) is homogeneous and can be integrated to derive
\be\label{v1}v(t)=(t-1)^{\beta}(t-B/A)^{-\beta},\ee where $\beta=
1/{(B-A)}$ is given by (\ref{1.6}).

Since in this case $\beta>1$, so by (\ref{4.7}), the equation
\be\label{v2}v(t)\phi_1(t)=A(t-1)^{1+\beta}(t-B/A)^{1-\beta}=0\ee
has a unique root $1$. Thus, the parameter $b=1$ satisfies
(\ref{4.6}).

Substituting the parameter $b$ and the function $v(t)$ into
(\ref{4.2}) and removing a constant multiplicative factor, we obtain
a solution $u_1(k)$ to (\ref{4.1}) for $k\geq 1$:
\be\label{4.9}u_1(k)=\int_0^1t^{k-1}\left(\frac{1-t}{1-\zeta
t}\right)^{\beta}dt,\ee where $\zeta=A/B$.

The order of the function $u_1(k)$ with respect to $k$ is given by
the following Lemma.

\bl$^{\makebox{\rm [Lemma 6.1 in \cite{CFV}]}}$\label{l4} Let $k\geq
1.$ Then \be\label{4.9'}u_1(k)=(1+O(k^{-1})){\cal
D}_1k^{-(1+\beta)}\ee for ${\cal D}_1={\cal D}_1(\a,\a_1)$ a fixed
constant.\el

In Case of $\a_1>\a_c$, we have $B<A$, and equation (\ref{4.5}) has
the same solution as (\ref{v1}). In addition, under the conditions
(\ref{1.3}) and (\ref{1.4}), one further has $\beta<-1$, and then
the equation (\ref{v2}) has a unique root $\gamma:=B/A$ as given in
(\ref{1.6}). So we can take $b=\gamma$ to satisfy (\ref{4.6}). Thus
 $$u_2(k)=\int_0^{\gamma}t^{k-1}\left(\frac{\gamma-t}{1-t}\right)^{-\beta}dt=\gamma^{k-\beta}\int_0^1t^{k-1}\left(\frac{1-t}{1-\gamma t}\right)^{-\beta}dt$$
 is a solution to (\ref{4.1}) for $k\geq 1$.

By Lemma \ref{l4}, we have
 \be\label{4.12}u_2(k)=(1+O(k^{-1})){\cal D}_2\gamma^kk^{-1+\beta} \ee for
 some fixed constant ${\cal D}_2={\cal D}_2(\a,\a_1)$.

 Finally, we consider the case of $\a_1=\a_c$. In this case $A=B$ and the equation
 (\ref{4.5}) can be integrated to derive
 $$v(t)=e^{-\mu/(1-t)}$$ with $\mu=1/A$ given in
 (\ref{1.6}). With same argument as in cases 1) and 2), take $b=1$ and define
 $$u_c(k)=\int_0^1t^{k-1}e^{-\frac \mu {1-t}}dt,$$ then $u_c$ is a
 solution to (\ref{4.1}) for $k\geq 1$.

Crudely, \be\label{4.11}u_c(k)\leq \int_0^1t^{k-1}dt=1/k.\ee The
precious representation of $u_c(k)$ can be found in Remark \ref{r1}.

Note that in all the three cases, $u_1$, $u_2$ and $u_c$ do not
satisfy equation (\ref{4.1}) when $k=0$. In fact, as calculated in
\cite{CFV}, for $i=1,2$ or $c$, we always have
\be\label{4.10}2A_2u_i(2)+(A_1+B_1)u_i(1)=[\phi_1(t)v(t)]^b_a=-\phi_1(0)v(0)\not=0.\ee

Now, we are going to solve (\ref{3.3}). By Remark \ref{r5}, we only
need to construct a solution for (\ref{3.3}) which satisfies the
requirements of Lemma \ref{l3}. Actually, we will construct such a
solution based on the solution of (\ref{4.1}) given above.

Denote by $g$ the solution for (\ref{4.1}), i.e., $g=u_1,u_2$ or
$u_c$ in the three cases respectively.

For $m>1$, define $w_k=0$ for $k\geq m$, $w_{m-1}=-\a_1/[(m-1)A_0]$
and for $j=m-2,m-3,\ldots,1$, let $w_j$ be such that
$$A_2(j+2)w_{j+2}+(A_1(j+1)+B_1)w_{j+1}+A_0jw_j=0.$$ Then $w_k$
satisfies (\ref{3.3}) for $k\geq 1$. Therefore, any linear
combination of $g$ and $w$ is a solution of (\ref{3.3}) for $k\geq
1$.

Now, let $$D=-\Df{2A_2w_2+(A_1+B_1)w_1}{2A_2g(2)+(A_1+B_1)g(1)} \ \
\makebox{and }\ d=-\Df{A_2(Dg(1)+w_1)}{B_1}.$$ Note that $D$ and $d$
depend on $g=u_1,u_2$ and $n_c$ respectively. By (\ref{4.10}), $D$
is well-defined.

Define $$d_k=\left\{\ba{ll}0,& \makebox{if }\ k=-1\\
d,& \makebox{if }\ k=0\\
Dg(k)+w_k, & \makebox{otherwise}\ea\right..$$ It is straightforward
to check that $d_k$ given above is the solution of (\ref{3.3}), by
(\ref{4.9'}), (\ref{4.12}) and (\ref{4.11}), we know that $d_k$
satisfies the requirements of Lemma \ref{l3}.

For $m=1$, we can take
$$D=-\Df{\a_1}{2A_2g(2)+(A_1+B_1)g(1)}, \ \
\ d=-\Df{DA_2g(1)}{B_1}$$ and directly define
$$d_k=\left\{\ba{ll}0,& \makebox{if }\ k=-1\\
d,& \makebox{if }\ k=0\ \ \ \ .\\
Dg(k), & \makebox{otherwise}\ea\right.$$ Similarly, in this case
$d_k$ is also a solution to (\ref{3.3}) which satisfies the
condition of Lemma \ref{l3}.

{\it Proof of Theorem \ref{th1}}: By the construction of the
solution $d_k$ and Lemma \ref{l3}, the theorem follows immediately
by taking
$$C_i=\left\{\ba{ll}-\Df{(2A_2w_2+(A_1+B_1)w_1){\cal D}_i}{2A_2u_i(2)+(A_1+B_1)u_i(1)},& \makebox{ for }
m>1\\[5mm]
-\Df{\a_1{\cal D}_i}{2A_2u_i(2)+(A_1+B_1)u_i(1)},& \makebox{ for }
m=1\ea\right.\ \ \makebox{for}\ i=1,2,c,$$where ${\cal D}_1$ and
${\cal D}_2$ are given in (\ref{4.9'}) and (\ref{4.12}), ${\cal
D}_c=1$.\QED

\section*{Acknowledgements} This work was begun when one of us
(Xian-Yuan~Wu) was visiting Institute of Mathematics, Academia
Sinica. He is thankful to the probability group of IM-AS for
hospitality. The authors thank Colin Copper for answering their
questions on the establishing of the recurrence (\ref{3.1}).

\end{document}